\documentclass[11pt]{article}

\usepackage{amsmath,amssymb,amsfonts,amsbsy,mathrsfs}

\renewcommand{\subsection}{\subsubsection}

\begin{document}

\title{\bf Elementary symmetrization of inviscid two-fluid flow equations giving a number of instant results}

\author{{\bf Lizhi Ruan}\\
Hubei Key Laboratory of Mathematical Physics, \\ School of Mathematics and Statistics, \\
Central China Normal University, Wuhan 430079, P.R. China \\
E-mail: rlz@mail.ccnu.edu.cn
\and
{\bf Yuri Trakhinin}\\
Sobolev Institute of Mathematics, Koptyug av. 4, 630090 Novosibirsk, Russia\\
and\\
Novosibirsk State University, Pirogova str. 2, 630090 Novosibirsk, Russia\\
E-mail: trakhin@math.nsc.ru
}

\date{ }

\maketitle

\begin{abstract}
We consider two models of a compressible inviscid isentropic two-fluid flow. The first one describes the liquid-gas two-phase flow. The second one can describe the mixture of two fluids of different densities or the mixture of fluid and particles. Introducing an entropy-like function, we reduce the equations of both models to a symmetric form which looks like the compressible Euler equations written in the nonconservative form in terms of the pressure, the velocity and the entropy. Basing on existing results for the Euler equations, this gives a number of instant results for both models. In particular, we conclude that all compressive shock waves in these models exist locally in time. For the 2D case, we make the conclusion about the local-in-time existence of vortex sheets under a ``supersonic'' stability condition. In the sense of a much lower regularity requirement for the initial data, our result for 2D vortex sheets essentially improves a recent result for vortex sheets in the liquid-gas two-phase flow.
\end{abstract}

\section{Introduction}

\label{s1}

 Two-phase or multi-phase flows are concerned with  flows with two or more components and have a wide range of applications in nature, engineering, and biomedicine.
 Examples include sediment transport, geysers, volcanic eruptions, clouds, rain in natural and climate system;
 mixture of oil and natural gas in extraction tubes of oil exploitation, oil transportation, steam generators, cooling systems, mixture of hot water and vapor of water in cooling tubes of nuclear power stations in energy production; bubble columns,  aeration systems, tumor biology, anticancer therapies, developmental biology,  plant physiology in  chemical engineering, medical and genetic engineering, bioengineering, and so on. Multi-phase flow  is much more complicated than  single-phase flow due to the existence of a moving and deformable interface and its interactions with multi-phases \cite{Brennen05,Ishii06,Kolev05,Kolev05-2}.

 In this paper, we consider two models of a compressible inviscid isentropic two-fluid flow.  The first one describes the liquid-gas two-phase flow. The second one can describe the mixture of two fluids of different densities or the mixture of fluid and particles. We now present corresponding systems of equations and review main existing results related to them.

\subsection{Inviscid liquid-gas two-phase flow}
We first consider the equations of the compressible inviscid liquid-gas two-phase isentropic flow:
\begin{equation}\label{1}
\left\{
\begin{array}{l}
 \partial_tm  +{\rm div}\, (m u )=0,\\[6pt]
  \partial_tn  +{\rm div}\, (n u )=0,\\[6pt]
 \partial_t(n u ) +{\rm div}\,(nu\otimes u ) + {\nabla}P=0,
\end{array}
\right.
\end{equation}
where $m=\alpha_g\rho_g$ and $n=\alpha_l\rho_l$ are the gas mass and the liquid mass respectively, the non-negative functions $\alpha_g$  and $\alpha_l$ satisfying the condition  $\alpha_g+\alpha_l=1$ denote the gas and liquid volume fractions respectively,    $\rho_g$ and $\rho_l$ are the gas and liquid densities,  $u$ denotes the mixed velocity of the liquid and the gas, and $P$ is the common pressure for both
phases.

In fact, the model \eqref{1} is a simplified version of the following general two-phase flow model:
\begin{equation*}
\left\{
\begin{array}{l}
\partial_tm+{\rm div}\,  (m u_g)=0,\\[3mm]
\partial_tn+{\rm div}\,  (n u_l)=0,\\[3mm]
\partial_t(m u_g+n u_l)+{\rm div}\,  (m u_g\otimes u_g+n u_l\otimes u_l)+\nabla P(m,n)=0,
\end{array}\right.
\end{equation*}
where $m$ and $n$ denote the same as above, and  $u_g$ and $u_l$ denote the gas and liquid velocities respectively. Just as in many papers concerned with viscous liquid-gas two-phase flow, in order to help further understanding the behavior of a solution to the
general model,  we use two densities and only one fluid velocity to study the simplified two-phase flow model, i.e., we restrict ourselves to a flow regime where the liquid and gas velocities can be assumed to be equal: $u_g=u_l=u$. Since the liquid phase is much heavier than the gas phase, we also neglect the gas phase in the mixture momentum equations.

In general, the pressure law $P=P(m,n)$ is a nonlinear complicated function of the densities. We will follow the simplification made in \cite{Evje08,Evje11SIAMJMA} which guarantees that the model is thermo-mechanically consistent in the sense that one can get a basic energy estimate easily. We assume that the pressure $P$ is  a smooth function of $(m, n)$ defined on $(0,+\infty)\times(0,+\infty)$.
As mentioned in \cite{Evje08,Evje11SIAMJMA}, depending on the realistic physical background, three typically different kinds of pressure law are respectively expressed as follows.

The first one has the simple form
\begin{equation}
P(m,n)= (\gamma -1) (m+n)^{\gamma},
\label{2}
\end{equation}
with $\gamma >1$ being a constant. Another one is
\begin{equation}
P(m,n)=C\left(\frac{m}{\rho_l-n}\right)^\gamma ,
\label{es2}
\end{equation}
with constants $C>0$ and $\gamma>1$, and the liquid density $\rho_l= {\rm const}$. The third one is taken in a little bit more complex form
\begin{equation}
P(m,n)=C\left(-b(m,n)+\sqrt{b(m,n)^2+c(m,n)}\right),
\label{es3}
\end{equation}
where $C=a_l^2/2$ and $k_0=\rho_0-({p_0}/{a_l^2})$ are positive constants,
$$
\begin{array}{l}
\displaystyle b(m,n)=k_0-n-\left(\frac{a_g}{a_l}\right)^2m=k_0-n-a_0m,\\[3mm]
\displaystyle c(m,n)=4k_0\left(\frac{a_g}{a_l}\right)^2m=4k_0a_0m,
\end{array}
$$
the constants $a_g$ and $a_l$ are reference sonic speeds in gas and liquid respectively, and the constants $p_0$ and $\rho_0$ are reference pressure and density.

Viscous two-phase flows have been investigated extensively, in particular, the existence, uniqueness, regularity,
asymptotic behavior, decay rate estimates,  and blow-up phenomena of solutions for various one-dimensional and
multi-dimensional viscous two-phase flows have been studied recently. We remark that the exhaustive literature list is beyond the
scope of the paper, and here we just mention some related results in \cite{Du15,Hao12,Wang16,Wen12,YaoL10,YaoL11JDE} (see also references therein) for the multi-dimensional viscous two-phase flow model (the viscous counterpart of system \eqref{1}). The theory of the inviscid two-phase flow \eqref{1} is comparatively mathematically underdeveloped although there have been some numerical studies;  see
\cite{Baudin05,Baudin05-2,Brennen05,Evje07} and references therein.

\subsection{Inviscid bi-fluid and fluid-particle flows}
We now consider another model of the compressible inviscid isentropic two-fluid flow:
\begin{equation}\label{3}
\left\{
\begin{array}{l}
 \partial_tm  +{\rm div}\, (m u )=0,\\[6pt]
  \partial_tn  +{\rm div}\, (nu )=0,\\[6pt]
 \partial_t((m +n) u ) +{\rm div}\,((m +n)u\otimes u ) +{\nabla}P=0,
\end{array}
\right.
\end{equation}
where $m$ and $n$ are the densities of two fluids and $u$ stands for the velocity of bi-fluid. The pressure $P=P(m,n)$ depending on two densities (possibly in a complicated way) is considered in the following form:
\begin{equation}
P=P(m,n)= m^{\alpha} +An^{\gamma}
\label{4}
\end{equation}
with three constants $\alpha \geq 1$, $\gamma \geq 1$ and $A>0$. System \eqref{3} is the inviscid version of a ``simple'' but still realistic bi-fluid model derived from physical considerations in \cite{Ishii06} (see also \cite{Bresch18}).

The viscous version of system \eqref{3} has been investigated extensively. In fact, if $\alpha =1$, $\gamma >1$ and $A=1$, by the asymptotic analysis of a system of coupled kinetic and fluid equations, namely, the Vlasov-Fokker-Planck equation and the compressible Navier-Stokes equations, the authors of \cite{Mellet07,Mellet08} proved the singular limit convergence towards the viscous version of  \eqref{3}.
This kind of model is widely used in the modeling of fluid-particle interactions (see, e.g., \cite{Laurent04}). In this case, its inviscid version \eqref{3} is a two-phase macroscopic model describing the motion of the mixture of fluid and particles, where $n$ is the density of the fluid and $m$ is the density of particles in the mixture. The derivation of the corresponding viscous two-phase macroscopic model from the Vlasov-Fokker-Planck equation and the compressible Navier-Stokes equations was first addressed by Carrillo and Goudon in \cite{Carrillo06}. Recently the existence of a weak solution to the viscous two-phase macroscopic model was proved in \cite{Novotny18,Vasseur17}. At the same time, according to our knowledge, there were almost no mathematical studies of the inviscid equations \eqref{3}.

For the above two models \eqref{1} and \eqref{3} of a compressible inviscid isentropic two-fluid flow, by introducing an entropy-like function, we reduce the equations of both models to a symmetric form which looks like the compressible Euler equations written in the nonconservative form in terms of the pressure, the velocity and the entropy. Basing on existing results for the Euler equations, this gives a number of instant results for both models. In particular, we conclude that all compressive shock waves in these models exist locally in time. For the 2D case, we make the conclusion about the local-in-time existence of vortex sheets under a ``supersonic'' stability condition. In the sense of a much lower regularity requirement for the initial data, our result for 2D vortex sheets essentially improves the recent result in \cite{HWY} for vortex sheets in the liquid-gas two-phase flow.

\section{Elementary symmetrization of two-fluid flow equations}

\label{s2}

Let us first consider system \eqref{1}. As in \cite{RWWZ}, we assume that
\begin{equation}
m >0\quad\mbox{and}\quad n > 0.
\label{4''}
\end{equation}
Our crucial idea is really very simple: we introduce the entropy-like function
\[
S=\frac{m}{n}> 0.
\]
It indeed plays the role of entropy and we will call it the ``entropy'' because, as for the usual entropy, we get the equation
\[
\frac{{\rm d} S}{{\rm d}t}=0
\]
with ${\rm d} /{\rm d} t =\partial_t+({u} \cdot{\nabla} )$, following from the first two equations of system \eqref{1}. Moreover, we can equivalently rewrite the first two equations of  \eqref{1} as
\begin{equation}
\frac{1}{\rho}\,\frac{{\rm d} \rho}{{\rm d}t} +{\rm div}\,u =0,\quad \frac{{\rm d} S}{{\rm d}t}=0,
\label{5}
\end{equation}
where for further convenience we use the notation $\rho :=n$.

We now recalculate the pressure $P$ in terms of $R$ and $S$. For the pressure law \eqref{2} we have
\begin{equation}
P=P(\rho ,S)= (\gamma -1) \rho^{\gamma} (S+1)^{\gamma}.
\label{8'}
\end{equation}
For the equation of state \eqref{es2} we obtain
\begin{equation}
P(\rho ,S)=C\left(\frac{\rho S}{\rho_l-\rho}\right)^\gamma ,
\label{8/2}
\end{equation}
whereas we rewrite \eqref{es3} as
\begin{equation}
P(\rho ,S)=C\left(-b(\rho ,S)+\sqrt{b(\rho ,S)^2+4k_0a_0\rho S}\right),
\label{8/3}
\end{equation}
where $b(\rho ,S):= k_0-\rho-a_0\rho S$.

In view of \eqref{5}, we obtain
\begin{equation}
\frac{{\rm d} \rho}{{\rm d}t}= \frac{1}{P_{\rho}}\,\frac{{\rm d} P}{{\rm d}t}.
\label{6}
\end{equation}
Then, taking into account \eqref{5} and \eqref{6}, we easily symmetrize the system of conservation laws \eqref{1} by rewriting it in the nonconservative form
\begin{equation}
\frac{1}{\rho P_{\rho}}\,\frac{{\rm d} P}{{\rm d}t} +{\rm div}\,{u} =0,\quad
\rho\, \frac{{\rm d}u}{{\rm d}t}+{\nabla}P  =0 ,\quad \frac{{\rm d} S}{{\rm d} t} =0.
\label{8}
\end{equation}
Equations \eqref{8} have absolutely the same form as the compressible Euler equations written in the nonconservative form. Writing down \eqref{8} in a matrix form, we get the symmetric system
\begin{equation}
A_0(U )\partial_tU+\sum_{j=1}^3A_j(U )\partial_jU=0,
\label{9}
\end{equation}
which is {\it hyperbolic} if $A_0>0$, where
\[
U= (P,u,S),\quad A_0= {\rm diag} (1/(\rho P_{\rho}) ,\rho ,\rho, \rho,1),
\]
\[
A_1=\left( \begin{array}{ccccc} \frac{u_1}{\rho P_{\rho}}&1&0&0&0\\[3pt]
1&\rho u_1&0&0&0\\
0&0&\rho u_1&0&0\\ 0&0&0&\rho u_1&0\\
0&0&0&0&u_1\\
\end{array} \right) ,\quad
A_2=\left( \begin{array}{ccccc} \frac{u_2}{\rho P_{\rho}}&0&1&0&0\\[3pt]
0&\rho u_2&0&0&0\\
1&0&\rho u_2&0&0\\ 0&0&0&\rho u_2&0\\
0&0&0&0&u_2\\
\end{array} \right),
\]
\[
A_3=\left( \begin{array}{ccccc} \frac{u_3}{\rho P_{\rho}}&0&0&1&0\\[3pt]
0&\rho u_3&0&0&0\\
0&0&\rho u_3&0&0\\ 1&0&0&\rho u_3&0\\
0&0&0&0&u_3\\
\end{array} \right).
\]

For the pressure law \eqref{2}, $P_{\rho}=\gamma P/\rho$ and, hence, the hyperbolicity condition $A_0>0$ holds under assumptions \eqref{4''}. For the pressure law \eqref{es2}, we have
\begin{equation}
P_{\rho}=\frac{\gamma\rho_l C(\rho S)^{\gamma -1}}{(\rho_l -\rho )^{\gamma +1}},
\label{7/2}
\end{equation}
whereas for the equation of state \eqref{es3} we calculate:
\begin{equation}
P_{\rho}= C\left( 1+ a_0S+\frac{\rho (a_0S+1)^2+k_0(a_0S-1)}{\sqrt{b(\rho ,S)^2+4k_0a_0\rho S}}\right).
\label{7/3}
\end{equation}
One can see that under assumptions \eqref{4''} we have $P_{\rho} >0$ in \eqref{7/2} and \eqref{7/3} (we omit simple calculations for \eqref{7/3}). That is, system \eqref{9} is again hyperbolic.

The elementary symmetrization of conservations laws \eqref{1} giving the quasilinear symmetric hyperbolic system \eqref{8}/\eqref{9}, which looks like the symmetric system of the compressible Euler equations, implies a number of results {\it almost for nothing}. For example, we have the local-in-time existence and uniqueness theorem \cite{Kato,VKh} in Sobolev spaces $H^s$ for the Cauchy problem for system \eqref{8}, with $s\geq 3$.

We now consider system \eqref{3}. Let
\begin{equation}
m \geq 0\quad\mbox{and} \quad n >0.
\label{4'}
\end{equation}
We again introduce the ``entropy'' $S=m/n\geq 0$. Introducing also the total density
\[
\rho =m +n >0,
\]
it is easy to see that the first two equations of \eqref{3} are equivalently rewritten as \eqref{5}. Clearly, we again have \eqref{6} with
\begin{equation}
P=P(\rho ,S) =\left( \frac{\rho S}{S+1}\right)^{\alpha} + A\left( \frac{\rho}{S+1}\right)^{\gamma}
\label{eqst}
\end{equation}
and
\begin{equation}
P_{\rho}=\frac{\partial P(\rho ,S)}{\partial \rho}=\frac{\alpha}{\rho}\left( \frac{\rho S}{S+1}\right)^{\alpha} +
\frac{\gamma A}{\rho}\left( \frac{\rho}{S+1}\right)^{\gamma} .
\label{7}
\end{equation}
That is, system \eqref{3} is also equivalently written as the symmetric hyperbolic system \eqref{8}/\eqref{9} (under the hyperbolicity condition \eqref{4'}).

In the next sections, using the symmetric form \eqref{8} of equations \eqref{1} and \eqref{3} and basing on existing results for the Euler equations, we get local-in-time existence results for shock waves and vortex sheets in two-fluid models \eqref{1} and \eqref{3}. We leave for the reader to find another important results taking place for the compressible Euler equations which could be directly (or almost directly) carried over models  \eqref{1} and \eqref{3} by making use of our simple idea of the introduction of the entropy-like function $S=m/n$.

\section{Shock waves}
\label{s3}

We consider system \eqref{1} for $t\in [0,T]$ in the unbounded space domain $\mathbb{R}^3$ and suppose that $\Gamma (t)=\{ x_1-\varphi (t,x')=0\}$ is a smooth hypersurface in $[0,T]\times\mathbb{R}^3$, where
$x'=(x_2,x_3)$ are tangential coordinates. We assume that $\Gamma (t)$ is a surface of strong discontinuity for the conservation laws \eqref{1}, i.e., we are interested in solutions of \eqref{1} that are smooth on either side of $\Gamma (t)$. To be weak solutions of (\ref{1}) such piecewise smooth solutions should satisfy the Rankine-Hugoniot jump conditions
\begin{equation}
[\mathfrak{j}_0]=0,\quad[\mathfrak{j}]=0,\quad \mathfrak{j}\left[u_{N}\right] +
|N|^2[P]=0,\quad \mathfrak{j}[{u}_{\tau}]=0
\label{16}
\end{equation}
at each point of $\Gamma$, where $[g]=g^+|_{\Gamma}-g^-|_{\Gamma}$ denotes the jump of $g$, with $g^{\pm}:=g$ in the domains
\[
\Omega^{\pm}(t)=\{\pm (x_1- \varphi (t,x'))>0\},
\]
and
\[
\mathfrak{j}_0^{\pm}=m (u_{N}^{\pm}-\partial_t\varphi),\quad \mathfrak{j}^{\pm}=n (u_{N}^{\pm}-\partial_t\varphi),\quad u_{N}^{\pm}=({u}^{\pm} {\cdot}{N}),\quad {u}^{\pm}_{\tau}=(u^{\pm}_{\tau _1},u^{\pm}_{\tau _2}),\quad
u^{\pm}_{\tau _i}=({u}^{\pm} {\cdot}{\tau}_i),
\]
\[
{N}=(1,-\partial_2\varphi,-\partial_3\varphi ),\quad
{\tau}_1=(\partial_2\varphi,1,0),\quad {\tau}_2=(\partial_3\varphi,0,1),\quad  \mathfrak{j}:=\mathfrak{j}^{\pm}|_{\Gamma},\quad
\mathfrak{j}_0:=\mathfrak{j}_0^{\pm}|_{\Gamma}.
\]

We now consider shock waves. For them there is a nonzero mass transfer flux across the discontinuity surface: $u_{N}^{\pm}|_{\Gamma}\neq 0$, i.e., $\mathfrak{j}\neq 0$ (as well as $j_0\neq 0$ and $j_0+j\neq 0$). Then, the first two conditions in \eqref{16} imply
\begin{equation}
[S]=0.
\label{17}
\end{equation}
Recalling our notation $\rho :=n$, by virtue of \eqref{16}, \eqref{17}, we get the following jump conditions for shock waves:
\begin{equation}
[S]=0,\quad[\mathfrak{j}]=0,\quad \mathfrak{j}\left[u_{N}\right] +
|N|^2[P(\rho ,S)]=0,\quad [{u}_{\tau}]=0,
\label{16'}
\end{equation}
where $\mathfrak{j}^{\pm}=\rho (u_{N}^{\pm}-\partial_t\varphi)$. Since the ``entropy'' $S$ is continuous, for a fixed $S$  and without the first condition system \eqref{16'} coincides with the Rankine-Hugoniot conditions for {\it isentropic} Euler equations with the equation of state $P=P(\rho ):=P(\rho ,S)|_{S\equiv{\rm const}}$.

According to the results in \cite{M1,M2,Met} and their extension to hyperbolic symmetrizable systems with characteristics of variable multiplicities \cite{Kwon,MZ}, all uniformly stable shocks are structurally stable. Roughly speaking (we do not discuss regularity, compatibility conditions, etc.), this means that if the uniform Lopatinski condition holds at each point of the initial shock, then this shock exists locally in time. In other words, as soon as planar shock waves are uniformly stable according to the linear analysis with constant coefficients, we can make the conclusion about structural stability of corresponding nonplanar shocks. In this sense, the linear analysis with constant coefficients plays crucial role for shock waves.

The free boundary problem for shock waves is the problem for the systems
\begin{equation}
A_0(U^{\pm})\partial_tU^{\pm}+\sum_{j=1}^3A_j(U^{\pm} )\partial_jU^{\pm}=0\quad \mbox{in}\ \Omega^{\pm}(t)
\label{21}
\end{equation}
(cf. \eqref{9}) with the boundary conditions \eqref{16'} on $\Gamma (t)$ and corresponding initial data for $U^{\pm}$ and $\varphi$ at $t=0$. We can reduce this problem to that in the fixed domains $\mathbb{R}^3_{\pm}=\{\pm x_1>0,\quad x'\in\mathbb{R}^2\}$ by the simple change of variables $\tilde{x}_1=x_1-\varphi (t,x')$, $\tilde{x}'=x'$. Dropping tildes, we get a nonlinear initial boundary problem in the half-spaces $\mathbb{R}^3_{\pm}$ with the boundary conditions \eqref{16'} at $x_1=0$. Without loss of generality we can consider the unperturbed shock wave with the equation $x_1=0$. Then, the linearization of the nonlinear initial boundary value problem about its constant solution, with $\widehat{U}^{\pm}=(\widehat{P}^{\pm},\hat{u}^{\pm},\widehat{S}^{\pm})={\rm const}$ and $\hat{\varphi} =0$, gives a linear constant coefficients problem. For the perturbation of the ``entropy'' (we again denote it by $S$) we obtain the separate problem (cf. \eqref{8}, \eqref{17})
\begin{equation}
\partial_tS^{\pm} +(\hat{u}^{\pm}\cdot S) =f^{\pm}\quad \mbox{for}\ x\in\mathbb{R}_3^{\pm},\qquad [S]=g \quad \mbox{at}\ x_1=0,
\label{S}
\end{equation}
where, as usual, we introduce artificial source terms (for the above equations they are $f^{\pm}(t,x)$ and $g(t,x')$) to make the linear problem inhomogeneous.

For problem \eqref{S} we easily deduce the energy identity
\begin{equation}
I(t)-\int\limits_0^t\int\limits_{\mathbb{R}^2}
\left.\left( [\hat{u_1}](S^-)^2 +2\hat{u}^+_1S^-g+\hat{u}^+_1g^2 \right)\right|_{x_1=0}\,{\rm d}x'{\rm d}s=I(0)+2\sum_{\pm}\int\limits_0^t\int\limits_{\mathbb{R}^3_{\pm}}f^{\pm}S^{\pm}{\rm d}x{\rm d}s,
\label{27}
\end{equation}
where
\[
[\hat{u}_1]=\hat{u}_1^+-\hat{u}_1^-\quad\mbox{and}\quad
I(t)=\sum_{\pm}\|{S}^{\pm}(t)\|^2_{L^2(\mathbb{R}^3_{\pm})}.
\]
For {\it compressive} shocks ($[\rho ]>0$), it follows from $[\mathfrak{j}]=0$ that $[u_N]<0$. That is, for the unperturbed planar compressive shock we have $[\hat{u}_1]<0$. Then, for compressive shocks, by standard simple arguments, from \eqref{27} we get the a priori estimate without loss of derivatives
\[
\begin{split}
\sum_{\pm}\Big(\|{S}^{\pm}\|_{L^2([0,T]\times\mathbb{R}^3_{\pm})}+ & \|S^{\pm}_{|x_1=0}\|_{L^2([0,T]\times\mathbb{R}^2)}\Big) \\  & \leq C\left\{ \sum_{\pm}\left(\|{S}^{\pm}_{|t=0}\|_{L^2(\mathbb{R}^3_{\pm})}
+\|f^{\pm}\|_{L^2([0,T]\times\mathbb{R}^3_{\pm})}\right) +\|g\|_{L^2([0,T]\times\mathbb{R}^2)}\right\},
\end{split}
\]
where $C>0$ is a constant.

Having in hand the above a priori estimate, we can say that the linearized problem satisfies the uniform Lopatinski condition as soon as its subproblem for the perturbations of $P^{\pm}$ and $u^{\pm}$ does. This subproblem totally coincides with the linearized constant coefficients problem for shock waves in compressible isentropic gas dynamics. Moreover, for the pressure law \eqref{2}  we have the polytropic gas equation of state (cf. \eqref{8'})
\[
{P}(\rho )=c{\rho}^{\gamma},
\]
for the unperturbed constant solution, where $c=(\gamma -1)({S}+1))^{\gamma}>0$ is a constant for a fixed $S=\widehat{S}^{+}=\widehat{S}^-$ (recall that $[\widehat{S}]=0$). It was shown by Majda \cite{M1} that all isentropic {\it compressive} shock waves in a polytropic gas are uniformly stable. Their local-in-time existence was proved in \cite{M2} (see also \cite{Met} for some improvements of the results in \cite{M2}). Basing on these results, we can thus make the conclusion that all compressive shock waves in the compressible inviscid liquid-gas two-phase isentropic flow are structurally stable (for the pressure law \eqref{2}).

For a fixed constant $S=\widehat{S}^{+}=\widehat{S}^-$, the equation of state \eqref{8/2} reads
\[
P(\rho )=c\left(\frac{\rho }{\rho_l-\rho}\right)^\gamma ,
\]
where the constant $c =CS^{\gamma}>0$. The above equation of state $P(\rho )$ is a convex function. The result in \cite{M1} was obtained not only for the polytropic gas equation of state but also for a general case satisfying a stability condition and the subsonic requirement behind of the shock. In particular, according to this result, all isentropic {\it compressive} shock waves are uniformly stable provided that the equation of state is convex. Hence, we conclude that for the pressure law \eqref{es2} all compressive shock waves for model \eqref{1} are structurally stable. In fact, the same is true for the pressure law \eqref{es3} because the equation of state \eqref{8/3} is a convex function of $\rho$ (we omit simple calculations).

For the two-fluid model \eqref{3}, the Rankine-Hugoniot conditions read:
\begin{equation}
[\mathfrak{j}_1]=0,\quad[\mathfrak{j}_2]=0,\quad \mathfrak{j}\left[u_{N}\right] +|N|^2[P]=0,\quad \mathfrak{j}[{u}_{\tau}]=0,
\label{16.2}
\end{equation}
where
\[
\mathfrak{j}_1^{\pm}=m (u_{N}^{\pm}-\partial_t\varphi),\quad \mathfrak{j}_2^{\pm}=n (u_{N}^{\pm}-\partial_t\varphi),\quad
\mathfrak{j}^{\pm}=\rho^{\pm} (u_{N}^{\pm}-\partial_t\varphi),\quad \mathfrak{j}:=\mathfrak{j}^{\pm}|_{\Gamma},
\]
and we recall that $\rho =m+n$ is the total density. For shock waves ($\mathfrak{j}\neq 0$), we again have the ``entropy'' continuity \eqref{17}. The jump conditions \eqref{16.2} can be then equivalently rewritten as \eqref{16'}. Repeating the above arguments made for model \eqref{1}, we reduce the question about structural stability of shock waves for equations \eqref{3} to that for shock waves for the compressible isentropic Euler equations with the equation of state (cf. \eqref{eqst})
\[
{P}({\rho})=c_1{\rho}^{\alpha}+ c_2{\rho}^{\gamma},
\]
for the unperturbed constant solution, where $c_1=({S}/({S}+1))^{\alpha}\geq 0$ and $c_2=A(1/({S}+1))^{\gamma}>0$ are constants for a fixed $S=\widehat{S}^{+}=\widehat{S}^-$. The above equation of state $P(\rho )$ is a convex function (for $\alpha\geq 1$ and $\gamma >1$). That is, we can conclude that all compressive shock waves for the two-fluid model \eqref{3} are structurally stable.

\section{Vortex sheets}
\label{s4}

For vortex sheets there is no mass transfer flux across the discontinuity surface: $u_N^{\pm}|_{\Gamma}=\partial_t\varphi$. For this case,
both the Rankine-Hugoniot jump conditions \eqref{16} and \eqref{16.2} (for models \eqref{1} and \eqref{3} respectively) imply the boundary conditions
\begin{equation}
\partial_t\varphi = u_N^{\pm},\quad [P]=0,
\label{vs}
\end{equation}
which are the same as that for vortex sheets in gas dynamics (full or isentropic).

For model \eqref{1}, the free boundary value problem for vortex sheets is the problem for \eqref{8} in the domains $\Omega^{\pm}(t)$ with the boundary conditions \eqref{vs} and initial data for $U^{\pm}$ and $\varphi$. This free boundary problem totally coincides with that for vortex sheets for the nonisentropic Euler equations studied in \cite{MoTre,MoTreWang}. Note that in \cite{MoTre,MoTreWang} the Euler equations are written in form \eqref{8} for the polytropic gas equation of state (for which $\rho P_{\rho}=\gamma P$). Recall that for the equation of state \eqref{8'} we also have $\rho P_{\rho}=\gamma P$. Even formally, by introducing the ``entropy'' function $s=\ln (S+1)$, we can reduce \eqref{8'} to the polytropic gas equation of state $p=p(\rho , s)= A\rho^{\gamma}e^{s}$, with $A=\gamma -1$ (clearly, the last equation in \eqref{8} stays valid for $s=\ln (S+1)$). On the other hand, the fact that vortex sheets were formally considered in  \cite{MoTre,MoTreWang} for the polytropic gas equation of state is absolutely unessential because the results in \cite{MoTre,MoTreWang} stay valid for a general equation of state (as for 2D vortex sheets for the isentropic Euler equations whose linear and structural stability was proved in \cite{CS1,CS2} under the ``supersonic'' condition).

That is, the linear and nonlinear (structural) stability results obtained in \cite{MoTre,MoTreWang} for 2D vortex sheets can be {\it directly} carried over 2D vortex sheets for equations \eqref{1}. Referring to \cite{MoTreWang}, we conclude that a vortex sheet $x_1=\varphi (t,x_2)$ for the 2D case of \eqref{1} (with $x=(x_1,x_2)$ and $u=(u_1,u_2)$) is structurally stable under the ``supersonic'' condition
\begin{equation}
[u_2] > \big( c_+^{\frac{2}{3}}+ c_-^{\frac{2}{3}}\big)^{\frac{3}{2}},\quad [u_2]\neq \sqrt{2}\,(c_+ +c_-),
\label{supson}
\end{equation}
where $c_{\pm}= c(\rho^{\pm} ,s^{\pm})$ and $c=c(\rho ,s)=\sqrt{P_{\rho}}=\sqrt{\gamma P/\rho}$ (note that the second condition in \eqref{supson} is not physical and appeared in \cite{MoTre} and \cite{RWWZ} due to a requirement of the applied technique). The existence theorem in \cite{MoTreWang} was proved for a finite (not necessarily short) time but under the condition that the initial discontinuity is close to a rectilinear discontinuity corresponding to a constant background solution $\widehat{U}^{\pm}=(\widehat{P}^{\pm},0,\hat{u}_2^{\pm},s^{\pm})$, $\hat{\varphi}=0$, provided that this solution satisfies the ``supersonic'' condition \eqref{supson} (see \cite{MoTreWang} for the whole details: exact assumptions on the initial data, regularity of solutions, etc.). Clearly, instead of the assumption that the initial data are
close to a piecewise constant solution in the existence theorem we could assume that the time of existence is sufficiently short (see, e.g., such a local-in-time existence theorem for current-vortex sheets \cite{T09}).

We can rewrite the speeds $c_{\pm}$ in \eqref{supson} in terms of the original unknowns $m$ and $n$:
\[
c_{\pm}=c(m^{\pm},n^{\pm}),\quad c(m,n)=\sqrt{\frac{\gamma (\gamma -1)(m+n)^{\gamma}}{n}} =\sqrt{\left(1+\frac{m}{n}\right)P_n},
\]
where $P_n=\partial P (m,n)/\partial n$. Now $c_{\pm}$ are written exactly in the same form as in \cite{RWWZ} where the linear stability of 2D vortex sheets for equations \eqref{1} was studied. Introducing the entropy-like function $S=m/n$ and just referring to \cite{MoTre,MoTreWang}, we easily recover not only the linear stability result in \cite{RWWZ} but also a recent structural stability result in \cite{HWY}.

Moreover, referring to the existence theorem in \cite{MoTreWang}, we even essentially improve the nonlinear result in \cite{HWY}. Indeed, the existence of solutions $(U^{\pm},\varphi )\in H^{\mu -7}([0,T]\times\mathbb{R}_{\pm}^2)\times H^{\mu -6}([0,T]\times\mathbb{R})$ was proved in \cite{MoTreWang} for initial data $(U^{\pm}_0,\varphi_0 )\in H^{\mu +1/2}(\mathbb{R}_{\pm}^2)\times H^{\mu +1}(\mathbb{R})$, with an integer $\mu\geq 13$. At the same time, as for current-vortex sheets in \cite{T09}, the existence theorem in \cite{HWY} is proved for initial data and solutions belonging to the anisotropic weighted Sobolev spaces $H^k_*$ (see, e.g., \cite{OSY,Sec00}; one has Sobolev's embedding $H^k_*\hookrightarrow H^{[k/2]}$), instead of the usual Sobolev spaces $H^k$. More precisely, the existence of solutions $(U^{\pm},\varphi )\in H_*^{\mu -1}([0,T]\times\mathbb{R}_{\pm}^2)\times H^{\mu}([0,T]\times\mathbb{R})$ was proved in \cite{HWY} for initial data $(U^{\pm}_0,\varphi_0 )\in H^{2\mu +15}(\mathbb{R}_{\pm}^2)\times H^{2\mu +16}(\mathbb{R})$, with an integer $\mu\geq 15$.

For the equations of state \eqref{es2} and \eqref{es3} as well as for model \eqref{3} we still can use the results obtained in \cite{MoTre,MoTreWang} and make the conclusion about the structural stability of 2D vortex sheets under the ``supersonic'' condition \eqref{supson} with $c_{\pm}=\sqrt {P_{\rho}(\rho^{\pm} ,S^{\pm})}$, where the value $P_{\rho}$ is given by \eqref{7/2},  \eqref{7/3} and \eqref{7} respetively.

\bigskip

\section*{Acknowledgements}
The research of L.Z. Ruan was supported in part by the Natural Science Foundation of China $\#$11771169, $\#$11331005, $\#$11871236, Program for Changjiang Scholars and Innovative Research Team in University $\#$IRT17R46, and the Special Fund for Basic Scientific Research of Central Colleges $\#$CCNU18CXTD04.


\begin{thebibliography}{100}

\bibitem{Bresch18}
Bresch D., Desjardins B., Ghidaglia J.M., Grenier E., Hilliairet M.
{\em Multifluid models including compressible fluids. Handbook of Mathematical
Analysis in Mechanics of Viscous Fluids}, Eds. Giga Y. et
Novotny A. (2018), pp. 52.

\bibitem{Baudin05}
Baudin M., Berthon C., Coquel F., Masson R., Tran Q.H. A relaxation method for two-phase flow models with hydrodynamic closure law.
{\it Numer. Math.} {\bf 99} (2005), 411--440.

\bibitem{Baudin05-2}
Baudin M., Coquel F., Tran Q.H. A semi-implicit relaxation scheme for modeling two-phase flow in a pipeline,
{\it SIAM J. Sci. Comput.} {\bf 27}(2005), 914--936.

\bibitem{Brennen05}
Brennen C. E. {\em Fundamentals of Multiphase Flow}, Cambridge University Press, New York, 2005.

\bibitem{Carrillo06}
Carrillo, J.A., Goudon, T. Stability and asymptotic analysis of a fluid-particle interaction model. {\it Comm. Partial Differ. Eqs.}  {\bf 31} (2006), 1349--1379.

\bibitem{CS1}
Coulombel  J.-F., Secchi, P. The stability of compressible vortex sheets in two space dimensions. {\it Indiana Univ. Math. J.} {\bf 53}  (2004), 941--1012.

\bibitem{CS2}
Coulombel J.-F., Secchi P. Nonlinear compressible vortex sheets in two space dimensions. \textit{Ann. Sci. Ecole Norm. Sup.} \textbf{41} (2008), 85--139.

\bibitem{Du15}
Du L., Zhang Q. Blow up criterion of strong solution for 3D viscous liquid-gas
two-phase flow model with vacuum. {\it Phys. D} {\bf 309} (2015), 57--64.

\bibitem{Evje07} Evje S., Flatten T. On the wave structure of two-phase model.
{\it SIAM J. Appl. Math.} {\bf 67}, 487--511.

\bibitem{Evje08}
Evje S., Karlsen K. H. Global existence of weak solutions for a viscous two-phase model.
{\it J. Differential Equations} {\bf 245} (2008), 2660--2703.

\bibitem{Evje11SIAMJMA}
Evje S. Weak solutions for a gas-liquid model relevant for describing gas-kick in oil wells.
{\it SIAM J. Math. Anal.} {\bf 43} (2011), 1887--1922.

\bibitem{Hao12}
Hao C., Li H. Well-posedness for a multidimensional viscous liquid-gas two-phase flow model. {\it SIAM J. Math. Anal. } {\bf 44} (2012), 1304--1332.

\bibitem{HWY}
Huang F., Wang D., Yuan D. Nonlinear stability and existence of vortex sheets for inviscid liquid-gas two-phase flow. arXiv:1808.05905.

\bibitem{Ishii06}
Ishii M., Hibiki T. {\em Thermo-fluid Dynamics of Two-Phase Flow}, Springer-Verlag, New York, 2006.

\bibitem{Kato}
Kato T. The Cauchy problem for quasi-linear symmetric hyperbolic systems. {\it Arch. Ration. Mech. Anal.} {\bf 58} (1975), 181--205.

\bibitem{Kwon}
Kwon B. Structural conditions for full MHD equations. {\it Quar. Appl. Math.} {\bf 7}  (2009), 593--600.

\bibitem{Kolev05}
Kolev N.I. {\em Multiphase Flow Dynamics. Vol. 1. Fundamentals}, Springer-Verlag, Berlin, 2005.

\bibitem{Kolev05-2}
Kolev N.I. {\em Multiphase Flow Dynamics. Vol. 2. Thermal and Mechanical Interactions}, Springer-Verlag, Berlin, 2005.

\bibitem{Laurent04}
Laurent F., Massot M., Villedieu P. Eulerian multi-fluid modeling for the numerical simulation of
coalescence in polydisperse dense liquid sprays. {\it J. Comput. Phys.} {\bf 194} (2004), 505--543.

\bibitem{M1}
Majda A. The stability of multi-dimensional shock fronts. {\it Mem. Amer. Math. Soc.} {\bf 41}(275), (1983).

\bibitem{M2}
Majda A. The existence of multi-dimensional shock fronts. {\it Mem. Amer. Math. Soc.} {\bf 43}(281), (1983).

\bibitem{Mellet07}
Mellet A., Vasseur A. Global weak solutions for a Vlasov-Fokker-Planck/Navier-Stokes system of equations. {\it Math. Models Methods Appl. Sci.} {\bf 17} (2007), 1039--1063.

\bibitem{Mellet08}
Mellet A., Vasseur A. Asymptotic analysis for a Vlasov-Fokker-Planck/compressible Navier-Stokes system of equations. {\it Commun. Math. Phys.} {\bf 281} (2008), 573--596.

\bibitem{Met}
M\'etivier  G. Stability of multidimensional shocks. In: {\it Advances in the theory of shock waves}, Freist\"{u}hler H., Szepessy A. (eds.), Progr. Nonlinear Differential Equations Appl. Birkh\"auser, Boston, {\bf 47} (2001), 25--103.

\bibitem{MZ}
M\'etivier G., Zumbrun K. Hyperbolic boundary value problems for symmetric systems with variable multiplicities. {\it J. Differential Equations} {\bf 211}  (2005), 61--134.

\bibitem{MoTre}
Morando A., Trebeschi P. Two-dimensional vortex sheets for the nonisentropic Euler equations: linear stability. {\it J. Hyper. Diff. Eqs.} {\bf 5} (2008), 487--518.

\bibitem{MoTreWang}
Morando A., Trebeschi P., Wang T. Two-dimensional vortex sheets for the nonisentropic Euler equations: nonlinear stability. arXiv:1808.09290.

\bibitem{Novotny18}
Novotn\'y A., Pokorn\'y M. Weak solutions for some compressible multicomponent fluid models. arXiv:1802.00798v1.

\bibitem{OSY}
Ohno M., Shizuta Y., Yanagisawa T. The trace theorem on anisotropic Sobolev spaces.  {\it T\^{o}hoku Math. J.} {\bf 46} (1994), 393--401.

\bibitem{RWWZ}
Ruan L., Wang D., Weng S., Zhu C. Rectilinear vortex sheets of inviscid liquid-gas two-phase flow: linear stability. {\it Commun. Math. Sci.} {\bf 14} (2016), 735--776.

\bibitem{Sec00}
Secchi P. Some properties of anisotropic Sobolev spaces. {\it Arch. Math.} {\bf 75} (2000), 207--216.

\bibitem{T09}
Trakhinin Y. The existence of current-vortex sheets in ideal compressible magnetohydrodynamics. \textit{Arch. Ration. Mech. Anal.} \textbf{191} (2009), 245--310.

\bibitem{Vasseur17}
Vasseur A., Wen H., Yu C. Global weak solution to the viscous two-fluid model with finite energy. arXiv:1704.07354v2.

\bibitem{VKh}
Volpert A.I., Khudyaev S.I. On the Cauchy problem for composite systems of nonlinear differential equations. {\it Math. USSR-Sb.} {\bf 16} (1972), 517--544.

\bibitem{Wang16}
Wang W., Wang W. Large time behavior for the system of a viscous liquid-gas two-phase flow model in $\mathbb{R}^3$.
{\it J. Differential Equations} {\bf 261} (2016), 5561--5589.

\bibitem{Wen12}
Wen H., Yao L., Zhu C. A blow-up criterion of strong solution to a 3D viscous liquid-gas two-phase flow model with vacuum.
{\it J. Math. Pures Appl.}  {\bf 97} (2012), 204--229.

\bibitem{YaoL10}
Yao L., Zhang T., Zhu C. Existence of asymptotic behavior of global weak solutions to a 2D viscous liquid-gas
two-phase flow model. {\it SIAM J. Math. Anal.} {\bf 42} (2010), 1874--1897.

\bibitem{YaoL11JDE}
Yao L., Zhang T., Zhu C. A blow-up criterion for a 2D viscous liquid-gas two-phase flow model.
{\it J. Differential Equations} {\bf 250} (2011), 3362--3378.


\end{thebibliography}
\end{document}